\documentclass[12pt]{article}
\usepackage{amsmath,amsfonts,amsthm,amssymb}
\begin{document}
\newtheorem{lem}{Lemma}
\newtheorem{teo}{Theorem}
\newtheorem{def1}{Definition}

\pagestyle{plain}
\title{Weak Solutions of Stochastic Differential Equations 
over the Field of $p$-Adic Numbers}
\author{Hiroshi Kaneko and Anatoly N. Kochubei\footnotemark
}

\footnotetext[1]{Partially 
supported by the Ukrainian State Fund for Fundamental Research,
Grant 10.01/004}

\date{}
\maketitle

\begin{abstract}
Study of stochastic differential equations on the field  of 
$p$-adic numbers was initiated by the second author and has been developed
by the first author, who proved several results for the $p$-adic case,
similar to the theory of ordinary stochastic integral with respect to L\'evy processes on
the Euclidean spaces.
In this article, we present an improved definition 
of a stochastic integral on the field and prove the 
joint (time and space) continuity of the 
local time for $p$-adic stable processes. Then we use the method of 
random time change to obtain sufficient conditions for the existence of
a weak solution of a stochastic differential 
equation on the field, driven by the $p$-adic stable process, 
with a Borel measurable coefficient.
\end{abstract}

\section{Introduction}

Stochastic processes on the field $\mbox{\boldmath $Q$}_p$ of 
$p$-adic numbers have been studied in many papers, for example, by
Albeverio-Karwowski \cite{AK1, AK2}, Albeverio-Zhao \cite{AZ1, AZ2, AZ3},  
Evans \cite{EV}, Fig\`a-Talamanca et al. \cite{FG1,FG2},
Ismagilov \cite{Ism}, Kaneko \cite{KN1,KN2},  Karwowski--Vilela Mendes \cite{KM}, 
Kochubei \cite{KC0,KC1, KC2,KC3}, Varadarajan \cite{Var} and Yasuda 
\cite{Y1,Y2,Y3} (here we do not mention papers on processes with $p$-adic time
and other related subjects).
In particular, these authors constructed and studied wide
classes of Markov processes on $\mbox{\boldmath $Q$}_p$ (most of 
the results can be extended easily to more general local fields).
Their infinitesimal generators are usually hypersingular integral 
operators; the first and simplest example is Vladimirov's 
fractional differentiation operator $D^\alpha$ \cite{V,VVZ}
corresponding to the $p$-adic $\alpha$-stable process.

\vspace{1em}

Alike for the case of the Euclidean space, for the field 
$\mbox{\boldmath  $Q_p$}$ of $p$-adic numbers, stochastic processes and analysis
are tightly related. Indeed, we can find 
descriptions and facts in \cite{V, VVZ,KC2} which are transplanted for describing 
probablisitic concepts on $\mbox{\boldmath  $Q_p$}$.
In \cite{KC1}, the second author initiated the
theory of stochastic differential equations on $\mbox{\boldmath  
$Q_p$}$, and in \cite{KN1} the first author made several assertions on
stochastic integrals on $\mbox{\boldmath  $Q_p$}$ 
which are similar to ordinary It\^o calculus on the Euclidean space.
On the other hand, Yasuda \cite{Y1} developed potential-theoretic notions 
related to $p$-adic L\'evy processes, in particular the stable 
process. An analytic potential theory over $\mbox{\boldmath  $Q_p$}$ 
was developed by Haran \cite{Har}.

\vspace{1em}

In \cite{KN1,KC1} we considered only strong solutions of 
stochastic differential equations on $\mbox{\boldmath  
$Q_p$}$ driven by the $p$-adic $\alpha$-stable process. Note that 
in the conventional theory of stochastic differential equations 
(over $\boldsymbol  {R}$), the study of weak solutions of such equations 
began only recently \cite{Z1,Z2,EK, KN}.

\vspace{1em}

In this paper, we initiate a theory of weak solutions for the 
$p$-adic case. Since we need to discuss stochastic differential
equations without continuity of the coeffficient as
assumed in \cite{KN1}, we first provide a more general construction of 
stochastic integrals than the one in the first author's paper. 
Our improved definition of stochastic integrals admits predictable integrand satisfying a finiteness condition
on its moment.
Then, we follow the approach by Zanzotto \cite{Z1} based 
on the method of random time change. In order to use this 
technique in our situation, we prove the joint 
(with respect to the space and time variables) continuity of the local 
time for the $p$-adic $\alpha$-stable process, as well as some facts 
on the occupation time on mesurable sets from $\mbox{\boldmath  $Q_p$}$
with positive Haar measure. After that, we apply the theory of the 
L\'evy system of a Markov process on a general state space 
\cite{RW}. In fact, this is one of the first applications of the general 
theory of Markov processes outside the usual Euclidean realm.
In Section 4, we give a sufficient condition for the existence of 
a weak solution of the stochastic differential 
equation

$$
X(t) = x + \int_0^t b(X_{s-})d Z(s)
$$

\noindent
on $\mbox{\boldmath  $Q_p$}$, with respect to the $\alpha$-stable process 
$\{Z(t)\}_{t \geq 0}$, with a locally bounded coefficient $b$. The conditions 
are much less restritive than those guaranteeing the existence of 
strong solutions of such equations \cite{KN1,KC1}.
\vspace{1em}

\vspace{1em}

\vspace{1em}

\section{Stochastic integral}
\vspace{1em} 

In this section, we will establish 
the notion of stochastic integrals
of predictable stochastic process
with respect 
to the L\'evy process  $\{Z(t)\}_{t \geq 0}$ with $Z(0) = 0$
determined by the sequence 
$A= \{a(m)\}$ satisfying 
\vspace{1em} 

\begin{itemize}
\item[\rm{(1)}]\quad
$a(m) \geq a(m+1)$,
\end{itemize}

\vspace{1em} 

\begin{itemize}
\item[\rm{(2)}]\quad
$\displaystyle\lim_{m \to \infty} a(m) = 0, 
\displaystyle\lim_{m \to -\infty} a(m) > 0 \;\mbox{ or } =
\infty, $
\end{itemize}

\vspace{1em} 
\noindent
and
\vspace{1em} 

\begin{itemize}
\item[\rm{(3)}]\quad
$\displaystyle\sum_{m = -\infty}^{\infty} a(m)p^{\gamma m} <\infty$
\end{itemize}
(see \cite{AK2,KC2,Y1} for the basic notions regarding L\'evy 
processes on $\mbox{\boldmath  $Q_p$}$).
Then, we can choose a filtration $\{\mathcal F_{t}\}_{t \geq 0}$ satisfying
$\mathcal F_{t} \subset \sigma [Z(s) \mid s \leq t]$ for any $t$
so that $\mathcal F_{t}$ is independent of 
$\sigma [Z(s + t) - Z(t) \mid s > 0]$ for every $t \geq 0$.

\vspace{1em} 
In this section, the smallest $\sigma$-field on $[0,\infty)\times \Omega$  based on which all left-continuous
$\{\mathcal F_t\}$-adapted processes are measurable will be denoted by $\mathcal S$.
The \mbox{\boldmath  $Q_p$}-valued stochastic process defined on $[0,\infty) \times \Omega$
is said to be {\it predictable} if it is measurable with respect to 
the $\sigma$-field $\mathcal S$.
As standard notations in the theory of Markov process, the starting point of $\{Z_t\}_{t \geq 0}$ will be indicated
in the notations for the probability measure and the expectation as $P_0$ and $E_0$ respectively.
\vspace{1em} 
\vspace{1em} 
\begin{lem}\quad {\sl The linear space $\Phi$ of the  
\mbox{\boldmath  $Q_p$}-valued bounded 
stochastic processes defined on $[0,\infty) \times \Omega$
satisfying the following conditions } \rm{(i)}\mbox{ and } \rm{(ii)}{\rm:}
 \vspace{1em} 
 \end{lem}

\begin{itemize}
\item[\rm{(i)}] \quad $\Phi$ contains all bounded left-continuous $\{\mathcal F_t\}$-adapted processes,
\item[\rm{(ii)}] \quad for any sequence $\{\phi^{(n)}\}_{n=1}^{\infty} \subset\Phi$, 
$\lim_{n \to \infty} \phi^{(n)}(t,\omega) = \phi(t,\omega)$ and
$\Vert \phi^{(1)} (t,\omega)\Vert_p 
\leq \linebreak \Vert \phi^{(2)} (t,\omega)\Vert_p\leq \Vert \phi^{(3)} (t,\omega)\Vert_p \leq \cdots $ 
imply $\phi \in \Phi$,
\end{itemize}

\vspace{1em} 

\noindent
contains all bounded predictable processes.

\vspace{1em} 
{\it Proof} \quad Since any bounded predictable process is described as the limit
of bounded $\mathcal S$-measurable simple functions, it suffices to prove
$1_{E} \in \Phi$ for any $E \in \mathcal S$. For the family 
$\mathcal S' = \{E \subset [0,\infty) \times \Omega \mid 1_E  \in \Phi\}$,
we can verify that
\vspace{1em} 

(i)\quad $[0,\infty) \times \Omega \in \mathcal S'$,

\vspace{1em} 

(ii)\quad $E_1, E_2 \in \mathcal S'$ and  $E_1 \subset E_2$ imply
$E_2 \setminus E_1 \in \mathcal S'$,

\vspace{1em} 

(iii)\quad $E_n \in \mathcal S'$ and  $E_1 \subset E_2 \subset \cdots$ imply
$\cup E_n \in \mathcal S'$.

\vspace{1em} 

For any finite set $\{Y_1(t), \cdots , Y_k(t) \}$ of 
bounded left-continuous $\{\mathcal F_t\}$-adapted processes
and finite set of balls $\{B_1,\cdots, B_k \}$ in \mbox{\boldmath  $Q_p$}, we easily
see $\bigcap\{(t,\omega)\mid Y_i(t,\omega)\in B_i \}\in \mathcal S'$. This is because
$1_{B_i}(Y_{i}(t))$ is bounded  $\{\mathcal F_t\}$-adapted left continuous process
for any $i = 1,2 \cdots, k$.

\vspace{1em} 

Since the family $\mathcal C$ of the sets described as $\bigcap\{(t,\omega)\mid Y_i(t,\omega)\in B_i \}$
with finite set \linebreak $\{Y_1(t), \cdots , Y_k(t) \}$ of 
bounded left-continuous $\{\mathcal F_t\}$-adapted processes
and finite set of balls $\{B_1,\cdots, B_k \}$ in \mbox{\boldmath  $Q_p$}
is closed under finite intersection and $\sigma[\mathcal C] = \mathcal S$,
we can derive from Dynkin's theorem (e.g., Lemma 5.1 in \cite{IW}) that $\mathcal S \subset \mathcal S'$. 
\qquad\qquad\qquad\qquad\qquad\qquad\qquad\qquad \qquad\qquad\fbox{}

\vspace{1em} 

For arbitrarily fixed $\gamma \geq 1$ and  $T > 0$, we denote by $\mathcal L^\gamma$ the family of \mbox{\boldmath  $Q_p$}-valued
predictable  processes each $\{\phi(t)\}_{t \geq 0}$ of which is adapted to
$\{\mathcal F_t\}_{t \geq 0}$ and satisfies
$E_0[\int_0^T \Vert\phi(t,\omega)\Vert_p^\gamma dt]< \infty$.
Here, we introduce a subfamily  $\mathcal L^0$ of $\mathcal L^\gamma$
each element of which admits the following expression:

$$
\phi(t,\omega) = f_0(\omega)1_{\{0\}}(t) 
+ \sum_{i=1}^{n-1}f_{i}(\omega)1_{(t_i, t_{i+1}]}(t)
$$

\vspace{1em}  
\noindent
with some sequence $\{f_{i}\}$ of random variables satisfying
$\Vert f_{i} \Vert_{L^{\infty}(\Omega;P)}< \infty$ and some
$0 = t_0 < t_1 < \cdots < t_n =T$.

\vspace{1em} 
\vspace{1em} 
{\bf Proposition 1}  \quad  {\sl $\mathcal L^0$ is a dense subfamily of $\mathcal L^\gamma$
with respect to the norm }

$$\Vert \phi \Vert_{\mathcal L^{\gamma}}
=\Big(\int_0^TE_0[ \Vert\phi(t,\omega)\Vert_p^\gamma] dt\Big)^{1/\gamma}.$$

\vspace{1em} 
\vspace{1em} 
{\it Proof}  \quad For any $\phi \in \mathcal L^\gamma$, the sequence $\{\phi^{(M)}\}_{M=1}^{\infty}$
of stochastic processes defined by $\phi^{(M)}(t,\omega) = \phi(t,\omega)
\times 1_{B(0,p^M)}(\phi(t,\omega))$ satisfies $\lim_{M \to \infty}
\Vert \phi - \phi^{(M)}\Vert_{\mathcal L^{\gamma}} = 0$. 
Therefore, we may assume that $\phi(t,\omega)$ vanishes
outside $B(0,p^M)$.

\vspace{1em} 

We define 
$\Phi$ to be the set of
$\phi\in \mathcal L^{\gamma}$ such that $\Vert\phi(t,\omega)\Vert_p\leq p^M$
 for any $ (t ,\omega) $  and there exists  a sequence 
$\{\phi_n\}_{n=1}^{\infty} \subset \mathcal L^0 $
 satisfying 
$\lim_{n \to \infty}\Vert \phi - \phi_n \Vert_{\mathcal L^{\gamma}}= 0$.
Then, the space $\Phi$ is a linear space and 
for any sequence $\{Z^{(n)}\}_{n=1}^{\infty} \subset\Phi$, 
$\lim_{n \to \infty} Z^{(n)}(t,\omega) = Z(t,\omega)$ and
$\Vert Z^{(1)}(t,\omega) \Vert_p
\leq \Vert Z^{(2)}(t,\omega) \Vert_p \leq \Vert Z^{(3)}(t,\omega) \Vert _p \cdots $imply $\{Z(t)\}_{t \geq 0} \in \Phi$.

\vspace{1em} 
If $\{\phi (t)\}_{t \geq 0}$ is a left-continuous $\{\mathcal F_t\}$-adapted process satisfying
$\Vert\phi(t,\omega)\Vert_p\leq p^M$ for any $(t ,\omega)$, then

$$
\phi_n(t,\omega) = \phi(0,\omega)1_{\{0\}}(t) 
+ \sum_{i=0}^{2^n}\phi(iT/2^n,\omega)1_{(iT/2^n, (i+1)T/2^n]}(t)
$$

\noindent
is in $\mathcal L^0$ and satisfies $\lim_{n \to \infty}
\Vert \phi_n - \phi\Vert_{\mathcal L^{\gamma}} = 0$.
Accordingly, Lemma 1 shows that $\Phi$ contains all $\{\mathcal F_t\}$-adapted
predictable
processes enjoying $\Vert\phi(t,\omega)\Vert_p \leq p^M$ for any $(t ,\omega)$.
\qquad \qquad \qquad \qquad \fbox{}

\vspace{1em} 
\vspace{1em} 
For any element
$\phi \in \mathcal L^{0}$ given by
$
\phi(t,\omega) = f_0(\omega)1_{\{0\}}(t) 
+ \sum_{i=1}^{n-1}f_{i}(\omega)1_{(t_i, t_{i+1}]}(t)
$
with some sequence $\{f_{i}\}$ of random variables satisfying
$\sup \Vert f_{i} \Vert_{L^{\infty}(\Omega;P_0)}< \infty$ and some
$0 = t_0 < t_1 < \cdots < t_n =T$, the stochastic integral $\int_0^t \phi(s) dZ(s)$ 
with respect to 
$\{Z(t)\}_{t \geq 0}$ is defined by
\vspace{1em}  
 
$$
\int_0^t \phi(s) dZ(s) =
\sum_{i=1}^{n-1}f_i(Z(t_{i+1}\land t) - Z(t_i\land t)) 
\qquad \mbox{ for }\quad 0 \leq t \leq T.
$$
\vspace{1em} 

\noindent 
The stochastic processes 
$\{\int_0^t \phi(s) dZ(s)\}_{t \in [0,T]}$ can be
regarded as a right continuous $\{\mathcal F_{t}\}$-adapted
process.

\vspace{1em}  
\vspace{1em}  

As in \cite{KN1}, we can show the following $:$
\vspace{1em}  
\vspace{1em} 
\begin{lem}\quad  {\sl  There exists a positive constant $C_{A,\gamma}$
satisfying the following  conditions} {\rm:}
\vspace{1em}  

\begin{itemize}
\item[\rm{(i)}] \qquad  $E_0\big[ \Vert Z(t) \Vert_{\scriptstyle p}^{^{\scriptstyle \gamma}} \big] \leq C_{A,\gamma}t$
\quad {\sl for all $t \geq 0$}.
\vspace{1em} 

\item[\rm{(ii)}]\qquad 
$E_0\Bigg[\displaystyle\sup_{0 \leq t \leq u}\Big\Vert \int_0^t \phi(s) dZ(s)
\Big\Vert_{\scriptstyle p}^{^{\scriptstyle \gamma}} \Bigg] \leq 
C_{A,\gamma}\int_0^{u}
E_0\big[\Vert \phi(s) \Vert_{\scriptstyle p}^{^{\scriptstyle \gamma}}\big]ds$ \quad
{\sl for }\;$0 \leq u \leq T.$
\end{itemize}
\end{lem}

\vspace{1em}

For any element $\phi\in \mathcal L^{\gamma}$, there exists a sequence
$\{\phi_n\}_{n=1}^{\infty} \subset \mathcal L^{0}$, such \linebreak
that $\lim_{n \to \infty}\int_0^T
E_0[\Vert \phi_n(t) - \phi(t) \Vert_{\scriptstyle p}^{^{\scriptstyle \gamma}}]dt = 0$. 
Then one can derive from Proposition 1 that
\vspace{1em} 

$$
\displaystyle\lim_{n, m \to \infty} 
E_0\Big[\sup_{0 \leq t \leq T}
\Big\Vert \int_0^t \phi_n(s) dZ(s) - \int_0^t \phi_m(s) dZ(s)
\Big\Vert_{\scriptstyle p}^{\scriptstyle \gamma}\Big] = 0.
$$
\vspace{1em} 

\noindent
Therefore, the stochastic
integral $\int_0^t \phi(s) dZ(s)$ of $\{\phi(t)_{t \geq 0}\}$ with respect to $\{Z(t)\}_{t \geq 0}$ 
can be defined as a unique  \mbox{\boldmath  $Q_p$}-valued process $\{Y(t)\}_{t \geq 0}$ satisfying
\vspace{1em} 

$$\displaystyle\lim_{n \to \infty}  
E_0\Bigg[\sup_{0 \leq t \leq T}\Big\Vert Y(t) - \int_0^t \phi_n(s) dZ(s)
\Big\Vert_{\scriptstyle p}^{\scriptstyle \gamma}\Bigg] = 0.$$

\vspace{1em}  
\vspace{1em}  
Let us denote by $\mathcal D([0,T] \to \mbox{\boldmath  $Q_p$})$ 
the space of all right continuous sample paths on the time interval $[0,T]$ to \mbox{\boldmath  $Q_p$}
with left limit at every point.
Since we already know that $\int_0^{\cdot} \phi_n(s) dZ(s)$ 
is a $\mathcal D([0,T] \to \mbox{\boldmath  $Q_p$})$-valued stochastic process, we immediately
see that $\int_0^{\cdot} \phi(s) dZ(s)$ is a 
$\mathcal D([0,T] \to \mbox{\boldmath  $Q_p$})$-valued variable as well.

\vspace{1em}  
\vspace{1em} 

Again as in \cite{KN1}, we will obtain wider perspectives
of $p$-adic stochastic integral so that it covers the one with
respect to the $\alpha$-stable 
processes. Indeed, we can consider
a random walk corresponding to
a sequence $A = \{a(m)\}$ satisfying (1), (2) and
\vspace{1em} 

\noindent
\begin{itemize}
\item[\rm{(4)}]\qquad \qquad\qquad$
\displaystyle\sum_{m = -\infty}^{0} a(m)p^{\gamma m} <\infty
\quad \mbox{ for a given real number }\gamma \geq 1.$
\end{itemize}

\vspace{1em} 
\noindent
In this case, defining sequences $A(M) = \{a(M;m)\}$ 
with properties (1), (2) and (3) by
\vspace{1em} 

$$
a(M;m) =
\left\{\begin{array}{ll}
a(m)& \quad\; \mbox{ if }m < M, \\
0& \quad\; \mbox{ if }m \geq M,
\end{array}\right.
$$
\vspace{1em} 

\noindent
we can prove the following proposition.

\vspace{1em} 
\vspace{1em} 

{\bf Proposition 2} \quad {\sl For any L\'evy process $\{Z(t)\}_{t \geq 0}$
corresponding to $A$ with properties {\rm(1)},{\rm (2)} and {\rm(4)}, we have the 
following} {\rm:}
\vspace{1em} 

(i)\quad {\sl The \mbox{\boldmath  $Q_p$}-valued process $\{Z(M;t)\}_{t \geq 0}$ defined by
$Z(M;t) = \int_{\mbox{\boldmath  $Q_p$}}\int_0^t f_M(z) 
\Lambda(ds,dz)$ is a random walk corresponding to 
$A(M)$, where
$\Lambda$ stands for the Poisson random measure of $\{Z(t)\}_{t \geq 0}$ and $f_M$ denotes the
\mbox{\boldmath  $Q_p$}-valued function defined by }

$$
f_M(z) = 
\left\{\begin{array}{ll}
p^m z& \qquad\mbox{ if }\;\;\Vert z \Vert_p = p^{m + M} \mbox{ with some integer }m \geq 0,\\
z& \qquad\; \mbox{{\sl otherwise}}.
\end{array}\right.
$$

\vspace{1em} 
(ii)\quad {\sl There exists a sequence $\{\Omega(M;T)\}_{M = 0}^{\infty}$ of events satisfying
$\lim_{M \to \infty} P_0(\Omega(M;T))  = 1$ and}

$$
\qquad\displaystyle\int_0^t \phi(s)dZ(M;s) = \int_0^t \phi(s)dZ(M + k;s) 
$$
\vspace{1em} 

\noindent
{\sl for any } $t \in [0,T]  $ {\sl and }  $k = 1,2,\cdots$ a.s. {\sl on } $\Omega(M;T)$
{\sl for any} $\phi$ of $\mathcal L^{\gamma}$.
\vspace{1em} 
\vspace{1em} 

(iii)\quad {\sl For any sequence }
$\{\phi(M;t)\}_{M = 0}^{\infty} \subset
\mathcal L^{\gamma}$ {\sl enjoying }
$ \phi(M;t) = \phi(M + k;t) $
{\sl for any }$ t \in [0,T] $ {\sl and } $k = 1,2,\cdots$
a.s. {\sl on }\;$\Omega(M;T)$,\\
$$\displaystyle\int_0^t \phi(M;s)dZ(M;s) = \int_0^t \phi(M + k;s)dZ(M + k;s) $$
{\sl for any } $ t \in [0,T] $  {\sl and } $k = 1,2,\cdots$
 a.s. {\sl on } $\Omega(M;T).$

\vspace{1em} 
\vspace{1em} 

{\it Proof} \quad The proof  can be done as in the proof of
Proposition 2 in \cite{KN1}.\qquad\qquad\qquad\qquad\qquad \fbox{}

\vspace{1em} 

Accordingly, if a $\mbox{\boldmath  $Q_p$}$-valued
predictable process $\{\phi(t)\}_{t \geq 0}$ adapted to
$\{\mathcal F_t\}_{t \geq 0}$ admits predictable processes $\{\phi(M;t)\}_{t \geq 0}$ 
from $\mathcal L^{\gamma}$
satisfying $\phi(t) = \phi(M;t)$ on $\Omega(M;T)$ for any $t \in [0,T]$ and
$M = 1,2,\cdots$,
we can define the stochastic integral $\int_0^t \phi(t) dZ(t)$ with 
respect to a stochastic process $\{Z(t)\}_{t \geq 0}$ determined by $A$ satisfying
(1), (2) and (4).
Indeed, the stochastic integral is defined as a unique right continuous 
$\{\mathcal F_{t}\}$-adapted process $\{Y(t)\}_{t \geq 0}$ satisfying 

$$
Y(t) = \int_0^t \phi(M;s)dZ(M;s) \quad\mbox{ a.s. on } \Omega(M;T) \mbox{ for all } M = 1,2,\cdots.
$$

\vspace{1em} 
\vspace{1em} 
\section{Local Time and Related Properties of $p$-Adic Stable Processes}

Let $\{Z(t)\}_{t \geq 0}$ be the symmetric $\alpha$-stable process on $\mbox{\boldmath  $Q_p$}$ with $\alpha >1$.
This stochastic process is characterized as the Hunt process  $\{Z(t)\}_{t \geq 0}$ satisfying $Z(0)=0$ with transition probability densities of the form
\begin{itemize}
\item[\rm{(5)}]$\displaystyle\qquad\qquad\qquad
P(t,x-y)=\int_{\mbox{\boldmath  $Q_p$}}\chi (-(x-y)\xi )e^{-t||\xi ||_p^\alpha}\mu(d\xi ),$
\end{itemize}
where $\chi$ is the canonical additive character on $\mbox{\boldmath  $Q_p$}$ and
$\mu$ is a Haar measure normalized by the requirement 
that the measure of a unit ball equals 1. The process $\{Z(t)\}_{t \geq 0}$ is given as 
the L\'evy process
determined by the sequence $A=\{a(m)\}$ of the form 
$$
a(m)=\frac{1-p^{-1}}{1-p^{-\alpha -1}}p^{-\alpha m}.
$$
Below $P_x$ will denote the probability law induced by $\{Z(t)+x\}_{t \geq 0}$ on
the trajectory space $\mathcal D([0,\infty) \to \mbox{\boldmath  $Q_p$})$ under $P_0$, and $ E_x$ will be the corresponding expectation. When the probability measure $P_x$ is provided,
the $\alpha$-stable process starting from $x$ will be denoted again by $\{Z(t)\}_{t \geq 0}$.

\vspace{1em} 

Here,  we can recall that Yasuda's result (\cite{Y1}) shows that $\alpha >1$ implies that every point $x\in \mbox{\boldmath  $Q_p$}$ is regular for $\{ x\}$. Accordingly, we have $ P_x\{ 
\tau_x=0\} =1$, where $\tau_x=\inf \{ t>0 \mid  Z(t)=x\}$ (see \cite{Y1}). Moreover, the 
Hunt process  $\{Z(t)\}_{t \geq 0}$ satisfies Hunt's conditions (A) and (F) in \cite{Hunt}, and
$ P_y\{ \tau_x<\infty\}=1$ for any $x,y \in \mbox{\boldmath  $Q_p$}$ (see \cite{Y1,KC2}).

\vspace{1em} 

Consider the random Borel measure 
$$
\nu (t,B)=\int_0^t I_B(Z(s))\,ds
$$
on $\mbox{\boldmath  $Q_p$}$, called the occupation time measure, where $I_B$ stands for the indicator of the Borel set 
$B$. If this random Borel measure admis a density $L_t^x$ with respect to the Haar measure $\mu$, 
that is 
$$
\int_BL_t^x\,\mu(dx)=\int_0^tI_B(Z(s))\,ds \quad a.s.
$$
for any Borel set $B$, then $\{L_t^x\}$ is called the local time of 
the process $\{Z(t)\}_{t \geq 0}$.

\vspace{1em} 

It is known \cite{BG1} that $\{Z(t)\}_{t \geq 0}$ in our situation admits the local time $\{L_t^x\}_{t \geq 0}$, for which
the function  $L_t^x$ is jointly measurable in $(t,x)$ and
continuous and monotone non-decreasing in $t$ for every $x$ with probability one. We note that 
$L_t^0>0$ is satisfied a.s. for any $t >0$ with probability one. Let us
prove almost sure joint continuity of $L_t^x$ in $(t,x)$.

\vspace{1em} 

Let $N,\delta$ be arbitrary positive numbers. It is known (see (V.3.28)
in \cite{BG2}) that for any $x,a,b\in \mbox{\boldmath  $Q_p$}$
\begin{itemize}
\item[\rm{(6)}]$\qquad\qquad\qquad\qquad\quad
P_x\left\{ \sup\limits_{0\le t\le N}\left| 
L_t^a-L_t^b\right| >2\delta \right\} \le 2e^Ne^{-\delta /\gamma_{a,b} },$
\end{itemize}
where $\gamma_{a,b} =[1-\psi_a(b)\psi_b(a)]^{1/2}$,
$\psi_a(x)=E_x\left( e^{-\tau_a}\right)$ and
$\psi_b(x)=E_x\left( e^{-\tau_b}\right)$.

\vspace{1em} 

For the function
$$
g^\lambda (x)=\int_0^\infty e^{-\lambda t}P(t,x)\,dt,$$

\noindent
with positive parameter $\lambda$, we can derive from (5) that
$$
g^\lambda (x)=\int_{\mbox{\boldmath  $Q_p$}}\frac{\chi (-x\xi 
)}{\lambda +||\xi||_p^\alpha }\mu(d\xi ).
$$
Since $\alpha >1$, the function $g^\lambda (x)$ is bounded and 
continuous in $x$, and
$ E_0\left( e^{-\tau_a}\right) =g^\lambda 
(a)/g^\lambda (0)$
(see Lemma 3.2 in \cite{Y1}).

\vspace{1em} 

Due to the spatial homogeneity, we have
$
\psi_a(b)= E_0\left( e^{-\tau_{a-b}}\right) =g^1 
(a-b)/g^1(0)\quad \mbox{ and }\quad \psi_b(a)=g^1(b-a)/g^1(0).
$
In our case, $g^\lambda (x)=g^\lambda (-x)$ is satisfied ($g^\lambda (x)$ depends
actually on $||x||_p$ as the Fourier transform of a radial function) and so we have
\begin{itemize}
\item[\rm{(7)}]\qquad\qquad\qquad\qquad$
\gamma_{a,b} =\left[ 1-\left( \frac{g^1(a-b)}{g^1(0)}\right)^2\right]^{1/2}
\le \sqrt{2}\left[ 1-\frac{g^1(a-b)}{g^1(0)}\right]^{1/2}.$
\end{itemize}

Next, we have the following lemma for the function
$$
h(||x||_p)=1-\frac{g^1(x)}{g^1(0)},\quad x\in \mbox{\boldmath  $Q_p$}.
$$

\vspace{1em} 
\begin{lem}
For any $N\in \mbox{\boldmath  $N$}$, there exists a positive constant $C_N$ such that
\begin{itemize}
\item[\rm{(8)}]\qquad\qquad\qquad$
0\le h(||x||_p)\le C_N||x||_p^{\alpha -1}\quad\mbox{ for any }x \mbox{ in } B(0, p^N). $
\end{itemize}
\end{lem}

\vspace{1em} 

{\it Proof.} We have 
$$
0\le h(||x||_p)\le C\int_{\mbox{\boldmath  $Q_p$}}\frac{1-\chi (-x\xi 
)}{1+||\xi ||_p^\alpha }\mu(d\xi )=C\int_{||\xi ||_p>||x||_p^{-1}}\frac{1-\chi (-x\xi 
)}{1+||\xi ||_p^\alpha }\mu(d\xi )=C\varphi (||x||_p)
$$
with some positive constant $C$, where 
$$\varphi (||x||_p)=\int_{||\xi ||_p>||x||_p^{-1}}\frac{1-\chi (-x\xi 
)}{1+||\xi ||_p^\alpha }\mu(d\xi ).
$$

For any $x \in  \mbox{\boldmath  $Q_p$}$ with $||x||_p=p^n$, from an integration formula in \cite{VVZ},
we can derive that
\begin{eqnarray*}
\varphi (||x||_p)&=&\sum\limits_{j=-n+1}^\infty \int_{||\xi ||_p=p^j}\frac{1-\chi (-x\xi 
)}{1+||\xi ||_p^\alpha }\mu(d\xi) \\
&=&\sum_{j=-n+1}^\infty 
\frac{1}{1+p^\alpha j}\left[ (1-p^{-1})p^j-\int_{||\xi ||_p=p^j}\chi (-x\xi 
)\,\mu(d\xi )\right] \\
&=&\frac{1}{1+p^{\alpha (1-n)}}\left[ (1-p^{-1})p^{-n+1}+p^{-n}\right] 
+(1-p^{-1})\sum\limits_{j=-n+2}^\infty \frac{p^j}{1+p^\alpha j}\\
&=&\frac{p^{-n+1}}{1+p^{\alpha (1-n)}}+(1-p^{-1})\sum\limits_{j=-n+2}^\infty \frac{p^j}{1+p^\alpha 
j}\\
&=&\frac{p}{p^n+p^\alpha \cdot p^{-n(\alpha -1)}}+(1-p^{-1})
\sum\limits_{k=2}^\infty \frac{p^{k-n}}{1+p^\alpha (k-n)}\\
&=&\frac{p}{p^n+p^\alpha \cdot p^{-n(\alpha -1)}}+(1-p^{-1})
\sum\limits_{k=2}^\infty \frac{p^k}{1+p^{\alpha k}\cdot p^{-n(\alpha 
-1)}}\\
&=&\frac{p}{||x||_p+p^\alpha ||x||_p^{-(\alpha -1)}}+(1-p^{-1})
\sum\limits_{k=2}^\infty \frac{p^k}{1+p^{\alpha k}||x||_p^{-(\alpha 
-1)}}\\
&=&||x||_p^{\alpha -1}\left\{ \frac{p}{||x||_p^{\alpha }+p^\alpha }+(1-p^{-1})
\sum\limits_{k=2}^\infty \frac{p^k}{||x||_p^{-(\alpha -1)}+p^{\alpha 
k}}\right\} .
\end{eqnarray*}

\vspace{1em} 
As $||x||_p\to \infty$, the expression in the braces tends to
$$
p^{1-\alpha }+(1-p^{-1})\sum\limits_{k=2}^\infty p^{k(1-\alpha )}=
p^{1-\alpha }+(1-p^{-1})\frac{p^{2(1-\alpha )}}{1-p^{1-\alpha }}.
$$
Therefore, we obtain the inequality (8) for any $x$ in $B(0,p^N)$.
\qquad\qquad\qquad\qquad\qquad\qquad\qquad \fbox{}

\vspace{1em} 
\begin{teo}
The function $(t,x)\mapsto L_t^x$ from $(0,\infty )\times \mbox{\boldmath  $Q_p$}$ to $(0,\infty )$ is continuous a.s. Moreover, for any 
$\kappa$ with $0<\kappa <(\alpha -1)/2$, any $T>0$ and $M\in 
\mbox{\boldmath  $Z_+$}$, there exists a random variable $C(\kappa 
,T,M)>0$ such that
\begin{itemize}
\item[\rm{(9)}]$\qquad\qquad\qquad\quad
\sup\limits_{0\le t\le T}\left| L_t^a-L_t^b\right| \le C(\kappa 
,T,M)||a-b||_p^\kappa \qquad a.s.$
\end{itemize}
for all $a,b$ in $B(0,p^M)$.
\end{teo}

{\it Proof.} Let us fix a positive integer $n$ and consider the set $S_n$ of
elements in $ \mbox{\boldmath  $Q_p$}$ of the form
$$
p^{-M}\left( \xi_0+\xi_1p+\cdots +\xi_{M+n}p^{M+n}\right),\quad \xi_j \in
\{ 0,1,\ldots ,p-1\}.
$$
Let us first prove that, for any $a,b\in S_n$ 
satisfying $||a-b||_p=p^{-n}$, Inequality (9) holds almost surely with some constant independent of $n$.
If $a,b\in S_n$ and 
$||a-b||_p=p^{-n}$, then it follows from (6), (7) and (8) that
$$
 P_0\left\{ \sup\limits_{0\le t\le T}\left| 
L_t^a-L_t^b\right| >2\delta_n\right\} 
\le 2e^N\exp{\Big(-\frac{\delta_n}{C||a-b||_p^{(\alpha -1)/2}}\Big)}
=2e^N\exp{(-\sigma p^{n((\alpha -1)/2-\kappa )})}\quad 
$$

\noindent
for any $\sigma >0$, where $\delta_n=p^{-\kappa n}$. 
\bigskip

The number of pairs $(a,b)$ of elements of $S_n$ satisfying
$||a-b||_p=p^n$ equals $p^{M+n+1}(p-1)/2$.
This is because the quantity of 
all elements of $S_n$ must be multiplied by the number $p-1$ of elements 
whose distance from a given element equals $p^{-n}$, and then divided 
by two, so that each pair to be counted only once. Since we have
$\rho_n :=
P_0\left\{ \sup\limits_{\substack{a,b\in S_n\\ 
||a-b||_p=p^{-n}}}\sup\limits_{0\le t\le T}\left| 
L_t^a-L_t^b\right| >2\delta_n\right\}\le (p-1)p^{M+n+1}e^N\exp{(-\sigma p^{n((\alpha -1)/2-\kappa 
)})}$,
we have
$\sum_{n=1}^\infty \rho_n<\infty$. By the 
Borel-Cantelli lemma, we see
$$
\sup\limits_{\substack{a,b\in S_n\\ 
||a-b||_p=p^{-n}}}\sup\limits_{0\le t\le T}\left| 
L_t^a-L_t^b\right| \le 2\delta_n,
$$
except a finite number of values of $n$ with probability one. In other words,
almost surely, there exists a positive integer $n_0$ such that

\begin{itemize}
\item[\rm{(10)}]$\qquad\qquad\qquad\qquad\qquad\qquad
\sup\limits_{0\le t\le T}\left| 
L_t^a-L_t^b\right| \le 2||a-b||_p^\kappa$
\end{itemize}
 
\noindent
for any $a$ and $b$ in $S_n$ satisfying $||a-b||_p=p^{-n}$
with some integer $n\ge n_0$.
\bigskip

Now, suppose that $a,b$ are arbitrary elements in $B(0,p^M)$
satisfying $||a-b||_p=p^{-n}$ with some integer $n\ge n_0$. Let us take 
canonical representations
$$
a=p^{-M}\left( \xi_0+\xi_1p+\cdots 
+\xi'_{M+n}p^{M+n}+\xi'_{M+n+1}p^{M+n+1}+\cdots \right) ,
$$
$$
b=p^{-M}\left( \xi_0+\xi_1p+\cdots 
+\xi''_{M+n}p^{M+n}+\xi''_{M+n+1}p^{M+n+1}+\cdots \right) 
$$

\noindent
for $a$ and $b$
with the coefficients from $\{ 0,1,\ldots ,p-1\}$ enjoying
$\xi'_{M+n}\ne \xi''_{M+n}$. For the sequences $\{a_k\}$ and $\{b_k\}$ given by
$$
a_k =
\left\{\begin{array}{ll}
p^{-M}\left( \xi_0+\xi_1p+\cdots +\xi'_{M+n}p^{M+n}\right) &\quad k = 0,\\
a_0+\sum\limits_{j=1}^k\xi'_{M+n+j}p^{n+j} &\quad k\ge 1,
\end{array}\right.
$$
\noindent
and
$$
b_k =
\left\{\begin{array}{ll}
p^{-M}\left( \xi_0+\xi_1p+\cdots +\xi''_{M+n}p^{M+n}\right) &\quad k = 0,\\
a_0+\sum\limits_{i=1}^k\xi''_{M+n+i}p^{n+i} &\quad k\ge 1,
\end{array}\right.
$$

\noindent
respectively,
we have
\begin{multline*} 
L_t^a-L_t^b=\left( L_t^a-L_t^{a_0}\right) +\left( 
L_t^{a_0}-L_t^{b_0}\right) +\left( L_t^{b_0}-L_t^b\right) \\
=\left( L_t^a-L_t^{a_1}\right) +\left( L_t^{a_1}-L_t^{a_0}\right) 
+\left( L_t^{a_0}-L_t^{b_0}\right) +\left( 
L_t^{b_0}-L_t^{b_1}\right) +\left( L_t^{b_1}-L_t^b\right) =\ldots 
\\
=\left( L_t^a-L_t^{a_k}\right) +\sum\limits_{j=1}^k 
\left( L_t^{a_j}-L_t^{a_{j-1}}\right) +\left( 
L_t^{a_0}-L_t^{b_0}\right) +\sum\limits_{i=1}^k 
\left( L_t^{b_{i-1}}-L_t^{b_i}\right) +\left( 
L_t^{b_k}-L_t^b\right) .
\end{multline*}

Since $a_k\to a$ and $b_k\to b$ as $k\to \infty$, we obtain that 
$L_t^a-L_t^{a_k}\to 0$ and $L_t^{b_k}-L_t^b\to 0$ in probability 
(see Corollary V.3.29 in \cite{BG2}). Passing to the limit, we find 
that 
$$
L_t^a-L_t^b=\left( L_t^{a_0}-L_t^{b_0}\right) +\sum\limits_{j=1}^\infty 
\left( L_t^{a_j}-L_t^{a_{j-1}}\right) +\sum\limits_{i=1}^\infty 
\left( L_t^{b_{i-1}}-L_t^{b_i}\right) \quad a.s.
$$
Hence, by (10), we get the inequality
\begin{eqnarray*} 
\left| L_t^a-L_t^b\right| &\le& 2||a-b||_p^\kappa +2\sum\limits_{j=1}^\infty 
p^{-(n+j)\kappa } +2\sum\limits_{i=1}^\infty p^{-(n+i)\kappa 
}\\
&=&2||a-b||_p^\kappa +\frac{4p^{-(n+1)\kappa }}{1-p^{-\kappa 
}}=\left( 2+\frac{4p^{-\kappa }}{1-p^{-\kappa }}\right)||a-b||_p^\kappa \quad a.s.
\end{eqnarray*}
for $a,b \in B(0,p^M)$ satisfying $||a-b||_p\le p^{-n_0}$. This is 
equivalent to (9). It is clear that the inequality (9) implies the 
required joint continuity. \qquad\qquad\qquad\qquad\qquad\qquad\qquad\qquad\qquad \fbox{}

\bigskip
Let us consider the behavior of the occupation time 
measure $\nu (t,B)$ as $t\to \infty$.

\medskip
\begin{teo}
For every Borel set $B\subset \mbox{\boldmath  $Q_p$}$ of a positive Haar 
measure,
\begin{itemize}
\item[\rm{(11)}]\qquad\qquad\qquad\qquad\quad
$ P_0\left\{ \lim\limits_{t\to \infty}\nu (t,B)=\infty 
\right\}  =1.$
\end{itemize}
\end{teo}

\medskip
{\it Proof}. Consider the transition probability semigroup $\{p_t\}$ given by
$$
(p_tf)(x)=\int_{\mbox{\boldmath  $Q_p$}}P(t,x-y)f(y)\,\mu(dy).
$$
It can be verified by a direct calculation that
$$
\int_{\mbox{\boldmath  $Q_p$}}(p_tf)(x)\,\mu(dx)=\int_{\mbox{\boldmath  $Q_p$}}f(x)\,\mu(dx)
$$
for the indicator function $f$ of any ball in $\mbox{\boldmath  $Q_p$}$. Thus this identity holds
for any locally constant function $f$ with a compact support, and therefore
for any $f\in L^2(\mbox{\boldmath  $Q_p$};\mu)$. This shows that the Haar measure $\mu$
is an invariant measure for the $\alpha$-stable process $\{Z(t)\}_{t \geq 0}$.

\vspace{1em} 
It is proved in \cite{Y1} (see also \cite{KC2}) that the $\alpha$-stable process $\{Z(t)\}_{t \geq 0}$ with $\alpha \geq 1$
is recurrent, so that it hits any open set in 
$\mbox{\boldmath  $Q_p$}$ after any large time (for these notions in the context of processes on general 
locally compact Abelian groups see \cite{PS}). Together with the 
existence of an invariant measure, this implies the Harris 
recurrence of $\{Z(t)\}_{t \geq 0}$ (Proposition X.3.9 in \cite{RY}), so that every 
Borel set $B$ of a positive Haar measure is recurrent. Due to
Propositions X.3.11 and X.2.2 in \cite{RY}, for the proof of (11),  it suffices to show that the functional
$$
M_B(f)=\lim\limits_{\beta \to \infty}\beta \int_{\mbox{\boldmath  $Q_p$}}
\mu(dx) \  E_x\int_0^\infty e^{-\beta 
t}f(x+Z(t))I_B(x+Z(t))\,dt,
$$
defined on bounded positive Borel measurable functions $f$ on 
$\mbox{\boldmath  $Q_p$}$, does not vanish. This follows from 
$$
M_B(f)=\lim\limits_{\beta \to \infty}\beta \int_{\mbox{\boldmath  $Q_p$}}
\mu(dx) \int_0^\infty e^{-\beta 
t}p_t(I_Bf)(x)\,dt=\int_Bf(x)\,\mu(dx)
$$
and $\mu(B)>0$.
\qquad\qquad\qquad\qquad\qquad\qquad\qquad\qquad\qquad\qquad\qquad\qquad\qquad\qquad\qquad\qquad \quad\fbox{}

\bigskip
The above properties of the local time make it possible 
to prove an analogue of the Engelbert-Schmidt 
zero-one law for our situation.

\medskip
\begin{teo}
Let $f$ be a non-negative Borel function on $\mbox{\boldmath  $Q_p$}$. The 
following properties are equivalent {\rm:}
\begin{itemize}
\item[\rm{(i)}] \quad$ P_0\big\{ \int\limits_0^tf(Z(s))\,ds<\infty,
\mbox{ for all } t\ge 0\big\} >0$.
\item[\rm{(ii)}] \quad$ P_0\big\{ \int\limits_0^tf(Z(s))\,ds<\infty,
\mbox{ for all } t\ge 0\big\} =1$.
\item[\rm{(iii)}] \quad$f$ is locally integrable on $\mbox{\boldmath  $Q_p$}$.
\end{itemize}
\end{teo}

\vspace{1em} 
\vspace{1em}

\section{Time change and moment estimate}

\vspace{1em} 

In this section, we will discuss a weak 
solution of 
the stochastic differential equation 
\vspace{0.5em} 

\begin{itemize}
\item[\rm{(12)}]\qquad $\displaystyle
\left\{\begin{array}{l}
d X(t) = b(X(t-))dZ(t ),\\
X(0 )= x,
\end{array}\right.$
\end{itemize}
\vspace{1em} 

\noindent
by taking some
$\alpha$-stable process on $\mbox{\boldmath  $Q_p$}$.
We will focus on the case that coefficient in the right-hand side
is given  by a $\mbox{\boldmath  $Q_p$}$-valued Borel measurable function $b$ defined
on $\mbox{\boldmath  $Q_p$}$.
\begin{def1} {\rm An $\{\mathcal F_t\}$-adapted  stochastic process $ \{ X(t)\}_{t \geq 0}$ defined on
a probability measure space $(\Omega, \mathcal F, P)$ is called a {\it solution of
the stochastic differential equation}
{\rm (12)}
if there exists an $\{\mathcal F_t\}$-adapted  $\alpha$-stable process $\{Z(t)\}_{t \geq 0}$ with $Z(0)=0$ satisfying }
\vspace{0.5em} 

$$
X(t) = x + \int_0^t b(X_{s-})d Z(s).
$$
\end{def1}
\vspace{1em} 
\begin{def1} {\rm A solution $\{X(t)\}_{t \geq 0}$ of equation {\rm (12)} is said to be {\it trivial}, if }
\vspace{0.5em} 

$$P_0(X(t) = X(0) \quad \mbox{ for all } \quad t \geq 0) = 1.$$

\end{def1}
\vspace{1em}

\begin{def1} {\rm For $x \in \mbox{\boldmath  $Q_p$}$, we say the coefficient
$b$ in {\rm (12)} satisfies Condition {\rm (H)} with respect to $x$ if}
\vspace{0.5em} 
$$\int_0^t \Big(\int_{B(0,p^L)}\frac{1}{\Vert b(x + y)\Vert_p^\alpha}P(s,y)\mu(dy)\Big)ds < \infty$$

\vspace{1em} 
\noindent
{\rm
for every integer $L$.}
\end{def1}

\vspace{1em} 

Here, we introduce an increasing process $\{C(t)\}_{t \geq 0}$ associated with 
the $\alpha$-stable process $\{Z(t)\}_{t \geq 0}$ defined by

$$
C(t) = \int_0^t \frac{1}{\Vert b(x + Z(s))\Vert_p^{\alpha}}ds.
$$
\vspace{1em} 
\begin{lem} \quad If $\alpha \geq 1$, then the increasing process $\{C(t)\}_{t \geq 0} $ satisfies the following 
{\rm :}

\begin{itemize}

\item[{\rm (i)}] \quad $P_0(C(t)< \infty) = 1$ for every $t \geq  0$.

\item[{\rm (ii)}] \quad $P_0(\lim_{t \to \infty}C(t) = \infty) = 1$.

\end{itemize}

\end{lem}
\vspace{1em} 

{\bf Proof} \quad For any integer $M$, we define an $\{\mathcal F_t\}$-stopping time
$\sigma_M = \inf \{t > 0 \mid  \Vert Z(t) \Vert_p \geq p^M\}.$
Then we have

$$
C(t) = C({t \wedge \sigma_M})1_{\{t < \sigma_M\}} + C(t) 1_{\{t \geq \sigma_M\}}.
$$

\vspace{1em}

\noindent
>From the fact that $P(\sigma_M >0)=1$, we can derive

\begin{eqnarray*}
E_0[C(t \wedge \sigma_M)1_{\{t < \sigma_M\}}]&=&E_0\big[1_{\{t < \sigma_M\}}\int_0^{t\wedge\sigma_M}\frac{1}{\Vert b(x + Z(s))\Vert_p^{\alpha}}ds\big] \\
&=&E_0\big[1_{\{\sup_{0 \leq s \leq t}\Vert Z(s) \Vert_p < p^M\}}\int_0^{t\wedge\sigma_M}\frac{1}{\Vert b(x + Z(s))\Vert_p^{\alpha}}ds\big] \\
&\leq&\int_0^t\bigg(\int_{B(0,p^M)}\frac{1}{\Vert b(x + y)\Vert_p^{\alpha}}P(s,y)\mu(dy)\bigg)ds.
\end{eqnarray*}

\vspace{1em}

Condition (H) imposed on the coefficient implies the finiteness of the right-hand side.
Since $C({t \wedge \sigma_M})$ is equal to $C(t)$ on $\{t < \sigma_M\}$ and $\sigma_M \to \infty$ as $M \to \infty$,
the first assertion has been proved.

\vspace{1em}

Since $1/\Vert b(x+y)\Vert_p^{\alpha}> \varepsilon$ for all $y$ with some $\varepsilon>0$, we can take 
$B_{\varepsilon}=\{y \in \mbox{\boldmath  $Q_p$} 
\mid 1/\Vert b(x +y)\Vert_p^{\alpha}\geq \varepsilon\}$  so that $\mu(B_{\varepsilon})>0$ is satisfied.
Therefore, the second assertion follows from Theorem 2 and the estimate
$\int_0^t1/\Vert b(x + Z(s))\Vert_p^{\alpha}ds= \int_{\mbox{\boldmath  $Q_p$}}
L_t^y/\Vert b(x + y)\Vert_p^{\alpha}\mu(dy) \geq \varepsilon \nu(t,B_{\varepsilon})$.
\qquad\qquad\qquad\qquad\qquad\quad
\fbox{}

\vspace{1em}

\vspace{1em}

Here, we define $\tau_t = \inf\{s \geq 0 \mid C(s) > t\}$ for every $t \geq 0$. Since Lemma 4 shows that
$\Vert b(x + Z(s))\Vert_p > 0$ except Lebesgue measure zero set in $[0,\infty)$ for every $\omega \in \Omega$
and $x \in  \mbox{\boldmath  $Q_p$}$, 
we have 

$$\tau_t 
= \int_0^{C_{\tau_t}}\Vert b(x + Z({\tau_s}))\Vert_p^{\alpha}ds = 
\int_0^t \Vert b(x + Z({\tau_s}))\Vert_p^{\alpha}ds.$$

\vspace{1em}

Let us now introduce an $\{\mathcal H_t\}$-adapted time changed process $\{Y(t)\}_{t \geq 0}$ defined by
$Y(t )=Z(\tau_t)$, where $\{\mathcal H_t\}_{t \geq 0}$ stands for the filtration given by $\mathcal H_t= \mathcal F_{\tau_t}$ 
for every $t \geq 0$. The objective of this section is
demonstrating that $X(t )= x + Y(t)$ is a solution of the stochastic differential equation (12) by assuming $b$ is locally bounded.

\vspace{1em}

Since the ball $B(0,p^L)$ centered at the origin and with radius $p^L$ is given as a disjoint 
union $\bigcup_{i}B(a_i, p^{\ell})$ of balls
centered at $a_i \in B(0,p^L), i = 1,2, \cdot \cdot \cdot, p^{L - \ell}$, we can define a Markov process $\{ Z_{L,\ell}(t)\}_{t \geq 0}$
by

$$
Z_{L,\ell}(t) =
\left\{\begin{array}{ll}
a_i &\quad Z(t )\in B(a_i, p^{\ell}),\\
p^{-(L+1)} &\quad   Z(t) \notin B(0,p^L).
\end{array}\right.
$$

\vspace{1em}
By applying [28, VI-28,3],  we can obtain the L\'evy system of the Markov process $\{Z_{L,\ell}(t) \}_{t \geq 0}$. From 
this observation, we can derive that the L\'evy system of the $\alpha$-stable process is given by $H(t) = t$
and $ N(x,dy) = (K/\Vert y-x \Vert_p^{1+\alpha})\mu(dy)$.
This shows 
that, for any non-negative $\mathcal B_{\Delta}$ measurable function $f$ on $(\mbox{\boldmath  $Q_p$})_{\Delta}$ satisfying $f(0) =0$
and any initial probability measure $m$,

\vspace{1em}
$$\tilde A^{f}(t)= K\int_{[0,t]}ds 
\int_{(\mbox{\boldmath  $Q_p$})_{\Delta}}\frac{f(y-Z(s-))}{\Vert y-Z(s-)\Vert_p^{1+\alpha}}\mu(dy)$$

\vspace{1em}

\noindent
is the dual predictable projection of $A_f(t) =  \int_{[0,t]}\int_{(\mbox{\boldmath  $Q_p$})_{\Delta}}f(y)\lambda(\omega,ds,dy)
=\sum_{s\leq t} f(X(s) -X(s-))$ under $P_m$, where $(\mbox{\boldmath  $Q_p$})_{\Delta}$  stands for the one point compactification of $\mbox{\boldmath  $Q_p$}$ with 
 the topological Borel $\sigma$-field $\mathcal B_{\Delta}$ and $\lambda(\omega,dt,dy)$ stands for 
 the jump measure $\sum_{s > 0}1_{\{Z(s,\omega)-Z_(s-,\omega) \not= 0\}}\linebreak
\delta_{(s,Z(s,\omega)-Z(s-,\omega))}(dt,dy)$ of the $\alpha$-stable process $\{Z(t)\}_{t \geq 0}$ with filtration $\{\mathcal F_t\}_{t \geq 0}$.

\vspace{1em}

Consider a time changed jump measure

$$\tilde \lambda = \lambda(\omega,\cdot)\circ (\tilde C(t))^{-1},$$

\vspace{1em}

\noindent
with respect to $\tilde C(t )= (C(t),x)$ with $C(t) = \int_0^t1/\Vert b(x + Z(s))\Vert^{\alpha}ds $.
Then, another representation of $\{Y(t)\}_{t \geq 0}$ is obtained. Indeed, it is not difficult to see that $Y(t )= \int_{[0,t]}\int_{\mbox{\boldmath  $Q_p$}} x \tilde \lambda(ds,dx)$.

\vspace{1em}

Similarly to Lemma 2.20 in \cite{Z1}, we can show the following {\rm :}

\vspace{1em}
\begin{lem} \quad For almost all $\omega$ in $\Omega$, $N_{\omega} = \{s \in [0,\infty) \mid b(X(s,\omega)) = b(x+ Y(s,\omega))=0\}$
has Lebesgue measure zero.
\end{lem}

\vspace{1em}

On the other hand, by performing 
time change given by $\tilde C(t)$, it turns out that
the stochastic process determined by this jump measure $\tilde \lambda$ is characterized by the 
L\'evy system given by $H_t = \int_0^t\Vert b(X(s)))\Vert^{\alpha}ds $ and
$N(x,dy) = (K/\Vert y -x\Vert_p^{1+\alpha})\mu(dy)$. Therefore,  for any
non-negative $\mathcal B_{\Delta}$ measurable  function $f$ on $(\mbox{\boldmath  $Q_p$})_{\Delta}$
and any initial probabaility measure $m$,

\vspace{1em}

$$\tilde A_{H}^{f}(t)= K\int_{[0,t]} \Vert b(X(s))\Vert^{\alpha}ds 
\int_{(\mbox{\boldmath  $Q_p$})_{\Delta}}\frac{f(y-X(s-))}{\Vert y -X(s-)\Vert_p^{1+\alpha}}\mu(dy)$$

\vspace{1em}

\noindent
is the dual predictable projection of $A_{H}^f(t) = \sum_{s\leq H_t} f(X(s)-X(s-))$ under $P_m$
(see [28, VI-28,1]).
By introducing a random measure $\tilde \pi(\omega,ds,dy) 
= \Vert b(X(s,\omega))\Vert_p^{\alpha}(K/\Vert y \Vert_p^{1+\alpha} )ds\mu(dy)$,
we have $\tilde A_{H}^{f}(t)= \int_{[0,t]}\int_{(\mbox{\boldmath  $Q_p$})_{\Delta}}f(y)\tilde \pi(\omega,ds,dy)$.
Accordingly, $\tilde \pi$ is the compensator of $\tilde \lambda$ in the sense given in \cite{IW}.

\vspace{1em}

Let us denote the topological Borel $\sigma$-field of $[0,\infty)$ by $\mathcal B[0,\infty)$ and
the topological Borel $\sigma$-field of $\mbox{\boldmath  $Q_p$}$ by $\mathcal B$.
Now, we define a map $\beta_{\omega}$ from $( [0,\infty)\times \mbox{\boldmath  $Q_p$} ,\mathcal B[0,\infty) \times \mathcal B)$
to $( [0,\infty)\times(\mbox{\boldmath  $Q_p$})_{\Delta}, \mathcal B[0,\infty) \times \mathcal B_{\Delta})$ given by

$$\beta_{\omega}(s,y) = \Big(\omega,s,\frac{y}{b(X(s-,\omega))}\Big)$$

\noindent
for every $\omega \in \Omega$.

\vspace{1em}
\begin{lem}\quad The random measure defined by $\Pi(\omega,dt,dy) = \beta_{\omega} (\tilde \lambda(\omega,dt,dy))$
is characterized by the compensator (intensity) $K\rho_{\Delta}$,
where $\rho_{\Delta}$ stands for the measure on $(\mbox{\boldmath  $Q_p$})_\Delta$ 
given as $(1/\Vert y \Vert_p) \mu(dy)$ on $\mbox{\boldmath  $Q_p$}$ and vanishing at $\Delta$
\end{lem}

\vspace{1em}

{\bf Proof} \quad Thanks to Theorem 2.1(a) in \cite{KL}, it suffices to show that $\beta_{\omega}(\tilde \pi) = (K  /\Vert x \Vert^{1 + \alpha})ds\mu(dx) $.
Since $b(X(s-)) = b(X(s))$ is satisfied for almost all $s$ with probability one, we have

\begin{eqnarray*}
&&\int_{[0,t]}\int_{(\mbox{\boldmath  $Q_p$})_{\Delta}}f(y)\beta_{\omega}(\tilde \pi)(ds,dy) \\
&&= 
\int_{[0,t]}1_{\{b(X(s-)\not=0\}}\int_{(\mbox{\boldmath  $Q_p$})_{\Delta}}f\Big(\frac{y}{b(X(s-,\omega)}\Big)
\Vert b(X(s-),\omega)\Vert_p^{\alpha}\frac{K}{\Vert y \Vert_p^{1+\alpha}}\mu (dy)
\end{eqnarray*}

\vspace{1em}

\noindent
with probability one.
By performing changing of the variables $z =y/b(X(s-,\omega))$ in the integral, we obatin

\begin{eqnarray*}
\qquad&&\int_{[0,t]}\int_{(\mbox{\boldmath  $Q_p$})_{\Delta}}f(y)\beta_{\omega}(\tilde \pi)(ds,dy) \\
&&= 
\int_{[0,t]}1_{\{b(X(s-)\not=0\}}ds\int_{(\mbox{\boldmath  $Q_p$})_{\Delta}}f(z)
\frac{K}{\Vert z \Vert_p^{1+\alpha}}\mu(dz)\\
&&= 
\int_{[0,t]}\int_{(\mbox{\boldmath  $Q_p$})_{\Delta}}f(z)
\frac{K}{\Vert z \Vert_p^{1+\alpha}}ds\mu(dz). \qquad\qquad\qquad\qquad\qquad\qquad\qquad\qquad\qquad \qquad\qquad \fbox{}
\end{eqnarray*}

\vspace{1em}

Now, we define an $\{\mathcal H_t\}$-adapted stochastic process $\{Z^{\ast}(t)\}_{t \geq 0}$ given by

$$Z^{\ast}(t)= \int_{[0,t]}\int_{\mbox{\boldmath  $Q_p$}}y \Pi(ds,dy).$$

\noindent
Then it turns out that  $\{Z^{\ast}(t)\}_{t \geq 0}$ is an $\alpha$-stable process.
This is because the random measure $\Pi$ admits the compensator $(K/\Vert z \Vert_p^{1+\alpha})ds\mu(dz)$.

\vspace{1em}

For $\Omega_1(M;T)=\{ \tau_T < \sigma_M\}$ and 
$\Omega_2(M;T)=\{\sup_{0\leq s\leq T}\Vert Z^{\ast}(s)\Vert  \leq p^M\}$, we have
$\lim_{M\to \infty}\linebreak P_0(\Omega_1(M;T)) =\lim_{M\to \infty}P_0(\Omega_2(M;T)) =1$.
By modifying $\{b(X(t-))\mid 0\leq t \leq T\}$ and $\{Z^{\ast}(t)\mid 0\leq t \leq T\}$
outside $\Omega(M;T) =\Omega_1(M;T) \cap \Omega_2(M;T)$ 
as the procedure established in Section 2,
we can define the stochastic integral of $b(X(t-))$ 
with respect to the $\alpha$-stable process $\{Z^{\ast}(t)\}_{t \geq 0}$.
Similarly to Lemma 2.26 in \cite{Z1}, we have the following assertion
on the stochastic integral $M(t) = \int_{[0,t]}b(X(s-))dZ^{\ast}(s)$:

\vspace{1em}
\begin{lem}\quad $P_0(Y(t) = M(t )\mbox{ for all } t) = 1$.
\end{lem}

\vspace{1em}

Consequently, the following assertion on solutions of (12) is now concluded.

\vspace{1em}
\begin{teo} \quad If the coefficient $b$ is locally bounded and if
{\rm (H) } is satisfied with respect to $x$,
then the stochastic differential equation {\rm (12) } admits 
a non-trivial solution.
\end{teo}

\vspace{1 em}

Here, we note that the $\alpha$-stable process $\{Z(t)\}_{t \geq 0}$ on $\mbox{\boldmath  $Q_p$}$
has an explicit description of the density $P(t,x)$ of transition probability kernel.
Indeed, we firstly recall that

$$
P_0(\Vert Z(t) \Vert \leq p^m)= \frac{p-1}{p}\sum_{i = 0}^{\infty}p^{-i}\exp(-p^{-\alpha m}t),
$$

\noindent
is obtained by the expression of transition probability found by Yasuda (1996, [3]).
Then, by denoting the right-hand side by $P_m(t)$, we immediately see that
$P(t,y)$ is given as
$((p-1)/p^{m+1})(P_m(t) - P_{m-1}(t))$ on $B(0,p^m) \setminus B(0,p^{m-1})$.

\vspace{1em} 

We  have the following assertion which gives a sufficient condition for Condition (H).

\vspace{1em}
{\bf Proposition 3.}  For the $\alpha$-stable process with 
$\alpha \geq 1$ on $\mbox{\boldmath $Q_p$}$, the integrability 
condition

$$
\int_{B(0,1)} \frac{1}{\Vert b(x+y)\Vert_p^{\lambda}}\mu(dy) < \infty
$$
with some positive number $\lambda$ satisfying $ \lambda > \alpha(1 + \alpha)$ implies

$$
\int_0^T\Big\{\int_{B(0,1)} \frac{1}{\Vert b(x+y)\Vert_p^{\alpha}}
P(t,y)\mu(dy) \Big\}dt < \infty.
$$
\vspace{1em}

{\it Proof}\hspace{1cm} For any positive real number $q$, one sees that
\begin{eqnarray*}
\int_{B(0,1)} P(t,y)^q\mu(dy) 
&=& \sum_{m=-\infty}^{0}p^{m}(P_m(t) - P_{m-1}(t))^{q}\frac{1}{p^{mq}}
\leq \sum_{m=-\infty}^{0}p^{m}P_m(t)^{q}\frac{1}{p^{mq}}\\
&\leq& C_1\sum_{m=-\infty}^{0}p^{m}\Big(p^{mq}
+ \Big(\sum_{i=0}^{-m}p^{-i}\exp(-p^{-\alpha(m + i)}t)\Big)^q\Big)\frac{1}{p^{mq}}\\
&\leq& C_1\sum_{m=-\infty}^{0}p^{m}\Big(p^{mq}
+ \Big(\sum_{i=0}^{-m}p^{-i}\frac{1}{1 + p^{-\alpha(m + i)}t}\Big)^q\Big)\frac{1}{p^{mq}}\\
&\leq& C_1\sum_{m=-\infty}^{0}p^{m}\Big(p^{mq}
+ \Big(\sum_{i=0}^{-m}p^{-i}\frac{p^{\alpha(m + i)}}{p^{\alpha(m + i)} + t}\Big)^q\Big)\frac{1}{p^{mq}}\\
&\leq& C_1\sum_{m=-\infty}^{0}p^{m}\Big(p^{mq}
+ \Big(\sum_{i=0}^{-m}\frac{p^{\alpha m}}{p^{\alpha m} + t}\Big)^q\Big)\frac{1}{p^{mq}}\\
&\leq& C_2\sum_{m=-\infty}^{0}p^{m}\Big(p^{mq}
+ \Big(\frac{ \Vert m\Vert^{q} p^{mq}}{(p^{\alpha m} + t)^{q}}\Big)\Big)\frac{1}{p^{mq}}\\
&\leq &C_2\sum_{m=-\infty}^{0}p^{m}\Big(1
+ \Big(\frac{\Vert m\Vert^{q} }{(p^{\alpha m} + t)^{q}}\Big)\Big)\\
&\leq& C_2\sum_{m=-\infty}^{0}p^{m}\Big(1
+ \Big(\frac{1}{(p^{\alpha m} + t)^{q}}\Big)\Big)
\end{eqnarray*}
\vspace{1 em}

\noindent
with some positive constants $C_1$ and $C_2$,
where $\alpha \geq 1$.
Accordingly, it is not difficult to see that


\begin{eqnarray*}
\int_0^T\Big\{\int_{B(0,1)} P(t,y)^q\mu(dy) \Big\}^{1/q}dt
&\leq& C_3\int_0^T
\Big\{\sum_{m=-\infty}^{0}p^{m}\Big(1
+ \Big(\frac{1}{(p^{\alpha m} + t)^{q}}\Big)\Big)\Big\}^{1/q}dt\\
&\leq& C_3\Big(\int_0^T
\sum_{m=-\infty}^{0}p^{m}\Big(1
+ \Big(\frac{1}{(p^{\alpha m} + t)^{q}}\Big)\Big)dt\Big)^{1/q}\\
&\leq& C_3\Big(T
+ \sum_{m=-\infty}^{0}p^{m}\int_0^T\frac{1}{(p^{\alpha m} + t)^{q}}dt\Big)^{1/q}\\
&\leq& C_3\Big(T
+ \sum_{m=-\infty}^{0}p^{m(1 - \alpha(q -1))}\Big)^{1/q}\\
\end{eqnarray*}
\vspace{1 em}
\noindent
with some positive constant $C_3$,
where the right-hand side converges, when $(\alpha + 1)/\alpha> q$. 

Since $\lambda > \alpha(1+\alpha)$, a real number $q'$ is 
introduced by $q'=\frac{\lambda}{\alpha}$ so that $q'> \alpha + 1$ is satisfied. Then, it turns out that
\begin{eqnarray*}
&&\int_0^T\{\int_{B(0,1)} \frac{1}{\Vert b(x+y)\Vert_p^{\alpha}}P(t,y)\mu(dy) \}dt \\
&\leq& (\int_{B(0,1)} \frac{1}{\Vert b(x+y)\Vert_p^{\alpha q'}}\mu(dy))^{1/q'}
(\int_0^T\{\int_{B(0,1)} P(t,y)^q\mu(dy) \}^{1/q}dt) \\
&\leq& (\int_{B(0,1)} \frac{1}{\Vert b(x+y)\Vert_p^{\lambda}}\mu(dy))^{1/q'}
(\int_0^T\{\int_{B(0,1)} P(t,y)^q\mu(dy) \}^{1/q}dt) < \infty
\end{eqnarray*}

\vspace{1 em}

\noindent
with the positive number $q$ given by $q = q'/(q'-1)$.
This is because $\lambda > \alpha(1 + \alpha)$ implies $(\alpha + 1)/\alpha> q$. 
\hspace{34 em}\fbox{}

\vspace{1em}

\vspace{1em}

\noindent
Department of mathematics \\
Tokyo University of Science \\
26 Wakamiya, Shinjuku, Tokyo\\
162-0827 Japan

\vspace{1em}
\vspace{1em}

\noindent
Institute of Mathematics \\
National Academy of Sciences of Ukraine \\
Tereshchenkivska 3, Kiev, 
01601 Ukraine

\end{document}